\documentclass[12pt]{article}
\usepackage{amssymb}
\usepackage{vmargin}
\usepackage[official]{eurosym}
\setpapersize{A4}
\setmarginsrb{1in}{17mm}{1in}{17mm}{1in}{0mm}{10pt}{9mm}
\usepackage{graphicx}
\usepackage{amsmath}
\usepackage{color}
\usepackage{enumerate}
\usepackage{array}

\usepackage{times}
\newtheorem{thm}{Theorem}[section]
\newtheorem{cor}[thm]{Corollary}

\newtheorem{defn}[thm]{Definition}
\newtheorem{rem}[thm]{Remark}

\newtheorem{lemma}{Lemma}

\numberwithin{equation}{section}
\begin{document}
\title{Admissible pair of spaces for not correctly solvable linear differential equations}
\author{N. Chernyavskaya, L. Dorel, L. Shuster}
\maketitle

\begin{abstract}
We consider the differential equation
\begin{align}\label{ab}
-y'(x)+q(x)y(x)=f(x), \quad x \in \mathbb R,
\end{align}
where $f \in L_{p}(\mathbb R)$, $p\in [1,\infty)$, and  $0\leq q \in
L_{1}^{\rm loc}(\mathbb R)$, $\int\limits_{-\infty}^{0}q(t)\,dt=\int\limits_{0}^{\infty}q(t)\,dt=\infty,$
\begin{align*}
q_{0}(a)=\inf_{x\in \mathbb R}\int_{x-a}^{x+a}q(t)\,dt=0 \quad{\rm \  for ~ any  }\quad a\in (0,\infty).
 \end{align*}
Under these conditions, the equation ({\rm \ref{ab}}) is  not correctly solvable
in $L_{p}(\mathbb R)$  for any $p \in [1, \infty) $.
Let $q^{*}(x)$ be the Otelbaev-type average of the function $q(t), t\in \mathbb{R}$,
at the point $t=x$;  $\theta(x)$ be a continuous positive function for $x \in \mathbb R$, and
\begin{align*}
L_{p,\theta }(\mathbb R) = \{f\in L_{p}^{\rm loc}(\mathbb R):\,
\int_{-\infty}^{\infty}|\theta(x)f(x)|^{p}\,dx<\infty
\},
\end{align*}
\begin{align*}
\|f\|_{L_{p,\theta}(\mathbb R)}=\left(\int_{-\infty}^{\infty}|\theta(x)f(x)|^{p}\,dx\right)^{1/p}\
\end{align*}

We show that if there exists a constant $c\in [1, \infty)$, such that the inequality
$$c^{-1}q^{*}(x)\leq \theta(x)\leq cq^{*}(x)$$
holds  for  all $x \in \mathbb{R}$, then under some additional conditions for $q$
the pair of spaces $\{L_{p, \theta}(\mathbb R); L_{p}(\mathbb R)\}$ is admissible for the
equation ({\rm  \ref{ab}}).
\\


\end{abstract}

\section{Introduction}{}

In the present paper we consider the equation
\begin{equation}\label{1.1}
-y'(x)+q(x)y(x)=f(x), \quad x \in \mathbb R,
\end{equation}
where $f \in L_{p} \, (L_{p}\,:= L_{p}(\mathbb R))$, $p\in [1,\infty)$, and
\begin{equation}\label{1.2}
0\leq q \in L_{1}^{\rm loc}(\mathbb R).
\end{equation}
By a solution of (\ref{1.1})
we mean any absolutely continuous function $y(x)$ that satisfies (\ref{1.1}) almost everywhere on $\mathbb R$. By $\theta$  we denote a continuous positive function of $x \in \mathbb R$ and, we set
\begin{align}\label{1.3}
L_{p,\theta }(\mathbb R) = \{f\in L_{p}^{\rm loc}(\mathbb R):\,
\int_{-\infty}^{\infty}|\theta(x)f(x)|^{p}\,dx<\infty
\},
\end{align}
\begin{align*}
\|f\|_{L_{p,\theta}(R)}=\left(\int_{-\infty}^{\infty}|\theta(x)f(x)|^{p}\,dx\right)^{1/p}.\
\end{align*}
Furthermore, for $\theta \equiv 1, ~ x\in \mathbb R$,  we write  $L_{p},~ ||\cdot||_{p}$  in place of $L_{p,\theta}$ and $||\cdot||_{p, \theta}$, respectively.\\
\begin{defn}\label{d1.1}
\label{D1}We say that the pair of spaces  $\{L_{p,\theta}(\mathbb R);L_{p}(\mathbb R)\}$  {\rm (}hence
$\{L_{p,\theta}; L_{p}\} ${\rm ) } is admissible for equation {\rm (\ref{1.1})} if

 {\rm (i)}
 for any $f \in L_{p}(\mathbb R)$, the equation {\rm (\ref{1.1})}  has a  unique solution $y\in L_{p,\theta};$

{\rm (ii)}
 there exists an absolute positive constant $c(p)\in (0, \infty)$ such that the solution
of {\rm \ref{1.1}}, $y \in L_{p,\theta}$  satisfies the inequality
\begin{equation}\label{1.4} \|y\|_{p,\theta}\leq c(p) \|f|\|_{p}, \quad {~ for \ any}~\quad f\in L_{p}.\textsf{}
\end{equation}
\end{defn}
Moreover, if $\theta \equiv 1$ and conditions  (i)--(ii)  are satisfied, we say that equation (\ref{1.1})
is {\it correctly solvable\/} in  $L_{p}$.

The problem of finding  minimal requirements for $q(\cdot)$  under which the conditions (i)--(ii) are
satisfied is studied in the case $\theta \equiv 1$ in \cite{Lukachev, Chernyavskaya}. The main result
of these works can be summarized as follows.

\begin{thm}\label{t1.2}{\rm \cite{Lukachev}} Let $p \in [1, \infty)$. Equation {\rm (\ref{1.1})} is correctly
solvable in $L_{p}$ if and only if there is an $a\in (0, \infty)$ such that
\begin{eqnarray}\label{1.5}
q_{0}(a)=\inf_{x\in \mathbb R}\int_{x-a}^{x+a}q(t)\,dt  >0.
\end{eqnarray}
\end{thm}

Thus, if conversely
\begin{equation}\label{1.6}
q_{0}(a)=0\quad \text{for  all}\  a\in (0, \infty),
\end{equation}
then for any $p \in [1, \infty)$  equation  {\rm (\ref{1.1})} is not correctly solvable in $L_{p}$.
In this paper we continue to study the problems from \cite{Lukachev, Chernyavskaya}. Our aim is to find an analogue
of Theorem \ref{t1.2} under the condition  (\ref{1.6}). To formulate the question precisely, we need a new statement,
which is equivalent to Theorem \ref{t1.2} (and easily follows from it, see Section 4).  To this end, we introduce the auxiliary function $d(x), ~ x \in \mathbb R$ (see  \cite{Chernyavskaya}). Assume that in addition to (\ref{1.2}) the following
condition is satisfied:
\begin{equation}\label{1.7}
\int_{-\infty}^{\infty}q(t)\,dt=\infty.
\end{equation}%
Then, for a fixed $x \in \mathbb R$, define
\begin{equation}\label{1.8}
d(x)\overset{\rm def}{=}\inf_{d \geq 0}  \{d: \int_{x-d}^{x+d}q(\xi)\,d\xi=2 \} .
\end{equation}%
We note that functions of type $d(x)$, $x \in \mathbb R,$ have been introduced and systematically
used by M. Otelbaev (see \cite{Otelbaev}).

 It follows immediately from ({\rm \ref{1.8}}) that the function
  \begin{align}\label{1.9}
  q^{*}(x)=\frac{1}{d(x)}, \quad x \in \mathbb R,
  \end{align}
   can be interpreted  as the Steklov-type average $Q(x,h)$ of the function $q(t)$
at the point $t=x$ with particular averaging step $h=d(x)$:
\begin{equation}\notag
Q(x,h) \overset{\rm def}{=} \frac{1}{2h}\int_{x-h}^{x+h}q(t)\,dt, ~ h>0
\Longrightarrow
\end{equation}%
\begin{equation}\notag
Q(x,d(x)) = \frac{1}{2d(x)}\int_{x-d(x)}^{x+d(x)}q(\xi)\,d\xi=\frac{1}{d(x)}=q^{*}(x).
\end{equation}
We call the function $q^{*}(x)$, which plays a significant role in this work, the Otelbaev-type average of the function $q(x), ~ x \in \mathbb R.$

\begin{thm}\label{t1.3} {\rm (}\textrm{see Section 4}{\rm)}.
Let $p \in [1, \infty)$ and suppose   {\rm (\ref{1.2})} holds. Then  equation {\rm (\ref{1.1})} is
correctly solvable in $L_{p}$ if and only if $(\ref{1.7})$ holds and $d_{0}~< \infty $ {\rm(}\textrm{equivalently},  $q_{0}^{*}>0{\rm )}$.
Here
\begin{eqnarray}\label{1.10}
d_{0}=\sup_{x\in \mathbb R}d(x) ~~ \left(q_{0}^{*}=\inf_{x\in \mathbb R}q^{*}(x)\right)
\end{eqnarray}
\end{thm}

\begin{cor}\label{c1.4}{\rm (}\textrm{see Section 4}{\rm)}.
Assume the conditions {\rm (\ref{1.2})} and {\rm  (\ref{1.7})} are satisfied.
Then $d_{0}< \infty$ if and only if $q_{0}(a)>0$ for some $a\in (0, \infty)$.
\end{cor}
Turning back to our problem, we conclude (see Theorem  \ref{t1.3}) that the equation (\ref{1.1})
is not correctly solvable in $L_{p}$ in two typical cases:

{\rm 1)}
   {\rm (\ref{1.7})} holds, but $d_{0}=\infty$,

{\rm 2)}
 {\rm (\ref{1.7})} does not hold.

In the present paper we study the first case. The second case will be considered in a forthcoming paper.

We are now in a position to formulate our problem in  precise terms.
Consider the equation (\ref{1.1}), with a function $q$ that satisfies the conditions (\ref{1.2}), (\ref{1.6}), (\ref{1.7})   and  $p\in[1, \infty)$.
The task is to find  a positive continuous  function $\theta(x), ~  x\in \mathbb R$ such that the pair
$\{L_{p,\theta}; L_{p}\}$ is admissible for equation (\ref{1.1}).
The solution of this problem (under some additional requirements on $q$ ) is presented in this paper.

The paper is divided into Sections 2 through 5. In Section 2 we present
preliminaries, Section 3 lists all our results, in Section 4 all
proofs are presented, and, finally, Section 5 is devoted to examples.

\textbf{Acknowledgment.} The authors thank Professors  Z. S. Grinshpun, E. Liflyand, and  J. Schiff for their
interest in our work and useful discussions.

\section{Preliminaries}{}
In the sequel (without special mention), we assume that the conditions ({\rm \ref{1.2}}), ({\rm \ref{1.7}}) are satisfied.
The symbols $c, c(\cdot), c_{1},c_{2},... $ stand for absolute positive constants, the values of which are not
essential for the exposition, and can even change within a single chain of calculations. One more definition is in order.

\begin{defn}{}Let $\varphi(x)$ ~and $\psi(x)$, $x \in \mathbb R$,  be continuous positive functions defined on interval
$(a,b)$, $-\infty\leq a<b\leq \infty$. We say that  $\varphi$ ~and $\psi$  are {\it weakly equivalent on\/} $(a,b)$,
and write $\varphi(x) \asymp \psi(x)$,  $x\in (a, b)$, if there exists an absolute constant $c\ge1$ such that
\begin{align*}
c^{-1}\varphi(x)\leq \psi(x)\leq c\varphi(x)\qquad {\rm for \ all } \qquad x\in (a,b).
\end{align*}
\end{defn}

\begin{thm}\label{t2.1}{\rm \cite{CSh}}
The  function $d(x)$  is finite and continuous on
$\mathbb R$, and possesses the following properties:
\begin{enumerate}
\item[\rm 1.] it satisfies the inequality
 \begin{align}\label{2.1}
 |d(x+h)-d(x)|\leq |h|, \qquad  whenever \qquad |h|\leq d(x);
 \end{align}
\item[\rm 2.]
it is differentiable almost everywhere on $\mathbb R;$
\item[\rm 3.]  for any $ x\in \mathbb R$, $\varepsilon \in [0,1]$,  and
$t \in [x-\varepsilon d(x), x+\varepsilon d(x)]$,
\begin{align}\label{2.2}
(1-\varepsilon)d(x)\leq d(t)\leq(1+\varepsilon)d(x)<(1-\varepsilon)^{-1}d(x).
\end{align}
\end{enumerate}
\end{thm}
Note also that we have the following implications:
\begin{equation}\label{2.4}
\int_{-\infty}^{0}q(t)\,dt=\infty ~ \Longrightarrow ~
\lim_{x \rightarrow -\infty}\left(x+d(x)\right)=-\infty,
 \end{equation}                   and
 \begin{equation}\label{2.5}
\int_{0}^{\infty}q(t)\,dt=\infty ~ \Longrightarrow ~ \lim_{x \rightarrow \infty}\left(x-d(x)\right)=\infty \end{equation}

 \begin{defn}\label{d2.2}{\rm \cite{CSh}}
Let  $x \in \mathbb R$,   $\varkappa(t)$  be  a positive continuous function defined on $\mathbb R$,
and $\{x_{n}\}^{\infty}_{n=1}$  $ ($resp., $\{x_{n}\}_{n=-\infty}^{-1})$ be a sequence of points.
Consider the intervals
\begin{equation}\notag
\bigtriangleup_{n}=[\bigtriangleup_{n}^{-},  \bigtriangleup_{n}^{+}],  \quad \bigtriangleup_{n}^{\pm}=x_{n}\pm \varkappa(x_{n}),\quad  n \geq 1 \quad (n\leq -1).
\end{equation}
We say that  the sequence  of intervals $\{\bigtriangleup_{n}\}^{\infty}_{n=1}$ $($resp., $\{\bigtriangleup_{n}\}_{n=-\infty}^{-1})$ forms an
$R(x, \varkappa)$-covering of $[x, \infty) $ $($resp., $(-\infty, x])$ if the following conditions are satisfied:
\begin{enumerate}
\item[\rm 1.] $ \bigtriangleup_{n}^{+}=\bigtriangleup_{n+1}^{-}$, $ n\geq 1$   $(resp., \  \bigtriangleup_{n-1}^{+}=\bigtriangleup_{n}^{-}$,  $n \leq -1);$
\item[\rm 2.] $ \bigtriangleup_{1}^{-}=x$  $(resp.,\   \bigtriangleup_{-1}^{+}=x);$
\item[\rm 3.]  $\bigcup_{n\geq 1}\bigtriangleup_{n}=[x, \infty )$  $(resp.,\bigcup_{n\leq -1}\bigtriangleup_{n}=[-\infty, x) ) $.
\end{enumerate}
\end{defn}

\begin{thm}\label{t2.3} {\rm \cite{CSh}}. Suppose that the positive continuous function $\varkappa$ defined on $\mathbb R$ satisfies the condition
\begin{equation}\label{2.6}
\lim_{t\rightarrow \infty}(t-\varkappa(t))=\infty  \quad  \big(resp., \lim_{t\rightarrow-\infty}(t+\varkappa(t) =-\infty\big).
\end{equation}
Then for every $x\in\mathbb R$, there exists an $R(x,\varkappa)$-covering of $[x, \infty)$  $($resp.,   of $(-\infty, x] )$.
\end{thm}

\begin{cor}\label{c2.4}{\rm \cite{CSh}}
Suppose that conditions {\rm (\ref{2.5})} {\rm((\ref{2.4}))} are satisfied. Then for every $x\in\mathbb
R$ there exists an $R(x,d)$-covering of $[x, \infty)$  $($resp.,   of $(-\infty, x] )$.
\end{cor}

\begin{defn}\label{d2.5}{\rm \cite{Shuster}}
Let $q$ be a function such that for some $a \geq 1, ~ b > 0 $ and
$x_0 \geq 1 $, for all $|x|\geq x_0$, we have the inequalities
\begin{align}\label{2.6}
a^{-1}d(x) \leq d(t) \leq ad(x) \quad{\rm ~ for }\quad  ~ |t-x|\leq bd(x).
\end{align}
Then we say that the function $q$ belongs to the class $K(\gamma)$ and write $q\in K(\gamma)$, where
\begin{align}\label{2.7}
\gamma=\gamma(a,b)=a\exp(-b/a^{2})
\end{align}
\end{defn}

According to Definition \ref{d2.5} and Theorem 26 (see \cite{Otel}), we have the following:
\begin{thm}\label{t2.6}{\rm \cite{Shuster}}
$q\in \textit{K}(\gamma)$  for any $\gamma >1$.
\end{thm}

\begin{lemma}\label{l2.7}{\rm {\cite{Shuster}}}
Let $a\geq1,~ b>0$ and~ $\gamma=\gamma(a,b)\leq e^{-1}$. Then $ b\geq1$.
\end{lemma}

The problem of studying whether $q \in K(\gamma)$ for a given value of parameter $\gamma$ is considered in the following two theorems.

\begin{thm}\label{t2.7} {\rm \cite{Shuster}}
Suppose a function $q$ can be decomposed into a sum
 \begin{equation}\label{2.8}
 q(x)=q_{1}(x)+q_{2}(x),  \quad x\in \mathbb R,
 \end{equation}
 where the function $q_{1}(x)$ is positive on $\mathbb R$ and absolutely continuous along with
 its derivative, and $q_{2} \in L_{1}^{\rm loc}(\mathbb R).$ Suppose further that
 \begin{align*}
 \varkappa_{1}(x)\rightarrow 0  \quad{\rm and}  \quad \varkappa_{2}(x)\rightarrow 0  \quad {\rm as}\quad |x|\rightarrow \infty ,
 \end{align*}
 where
\begin{equation}\label{2.9}
\varkappa_{1}(x)=\frac{1}{q_{1}(x)^{2}}\sup_{|\xi|\leq 2/q_{1}(x)} \left\vert\int_{x-\xi}^{x+\xi}q_{1}''(s)\,ds\right\vert , \quad x \in \mathbb R
\end{equation}
and
\begin{equation}\label{2.10}
\varkappa_{2}(x)=\sup_{|\xi|\leq 2/q_{1}(x)} \left\vert\int_{x-\xi}^{x+\xi}q_{2}(s)\,ds\right\vert , \quad x \in \mathbb R.
\end{equation}
Then  the following relations hold:
\begin{equation}\label{2.11}
q_{1}(x)d(x)=1+\varepsilon(x), \qquad  |\varepsilon(x)|\leq\varkappa_{1}(x)+\varkappa_{2}(x), \qquad |x|\gg 1,
\end{equation}
\begin{equation}\label{2.12}
q^{*}(x) \asymp q_{1}(x),\quad    x \in \mathbb R.
 \end{equation}
\end{thm}

\begin{thm}\label{t2.8} {\rm \cite{Shuster}}
Let the hypotheses of Theorem {\rm \ref{t2.7}}  hold. Suppose in addition that
\begin{align}\label{2.13}
\lim_{|x|\rightarrow\infty}\frac{q_{1}'(x)}{q_{1}^2(x)}=0
\end{align}
\begin{align}\label{2.14}
\lim_{|x|\rightarrow\infty}|x|q_{1}(x)=\infty
\end{align}
Then for any $\gamma_{0}\in(0,1)$,
\begin{align}\label{2.15}
q\in K(\gamma) ,{\rm ~ whenere } ~ \gamma =\gamma(a,b)\leq \gamma_0.
\end{align}
\end{thm}

\section{Results}{}

The next statement is the main result of this work.
\begin{thm}\label{t3.1} Consider the equation ({\rm \ref{1.1}}). Suppose  the coefficient $q$ satisfies the conditions {\rm(\ref{1.2})} {(\rm \ref{1.6})}, and
\begin{align}\label{3.1}
\int\limits_{-\infty}^{0}q(t)\,dt=\int\limits_{0}^{\infty}q(t)\,dt=\infty,
\end{align}
and, in addition, $q\in K(\gamma),~ \gamma \leq e^{-1}$.
 Then, if the function $\theta$ satisfies the condition
 \begin{align}\label{3.2}
\theta(x) \asymp q^{*}(x), ~ ~ x\in \mathbb{R},
\end{align}
the pair of spaces $\{L_{p,\theta}; L_{p}\}$  is admissible for the
equation {\rm(\ref{1.1})} for any $p \in [1, \infty)$.
\end{thm}

\begin{rem}\label{r3.2} {\rm We will check below the hypotheses  of Theorem {\rm\ref{t3.1}} {\rm(}\textrm{in view of theorems {\rm\ref{t2.7}} and {\rm\ref{t2.8}}}{\rm)} for a particular equation {\rm (see Section 5)}.}
\end{rem}

Theorem \ref{t3.1} is obtained from  Theorem \ref{t3.3}, which is  of considerable significance of its own.
To formulate Theorem \ref{t3.3} we introduce the following notations:
\begin{equation} \label{3.3}
\textit{I}(x)=\int_{-\infty}^{x}\exp\bigg\{
-\int_{t}^{x}q(\xi)\,d\xi\bigg\}\, dt,\quad x \in \mathbb R,
\end{equation}
\begin{equation}\label{3.4}
\textit{J}(x)=\int_{x}^{\infty}\exp
\bigg\{-\int\limits_{x}^{t}q(\xi)\,d\xi\bigg\}\, dt,\quad x \in \mathbb R,
\end{equation}
\begin{equation}\notag
\textit{S}(x)=\int_{-\infty}^{\infty}\exp\bigg\{-\bigg\vert \int_{x}^{t}q(\xi)\,d\xi\bigg\vert\bigg\} \, dt,\quad x \in \mathbb R.
\end{equation}

\begin{thm}\label{t3.3}
Suppose conditions {(\rm\ref{1.2})}, {(\rm\ref{3.1})}, and  $q \in K(\gamma), ~ \gamma  \leq e^{-1}$ are satisfied.
Then
\begin{align}\label{3.5}
J(x)~\asymp~I(x)~\asymp~S(x)~\asymp~d(x)~=~\frac{1}{q^{*}(x)},~\qquad  x \in \mathbb R.
\end{align}
\end{thm}
\begin{rem}\label{r3.4}
{\rm For $d_{0} < \infty$ (\textrm{see} {\rm(\ref{1.10})}), the relations {\rm (\ref{3.5})} were  proven in {\rm \cite{Shuster}}.
It follows from Corollary {\rm \ref{c1.4}} that the inequality  $d_{0}<\infty$ is incompatible with the requirement {\rm (\ref{1.6})}. Thus, we need to derive {\rm  (\ref{3.5})} in such a way that $d_{0}<\infty$ is not used {\rm ( see Section 4)}.}
\end{rem}

Below we give an example of an application of ({\rm {\ref{3.5}}).
Consider the integral
\begin{equation}\label{3.6}
F(x)=\int_{-\infty}^{\infty} G(x,t)\,dt,\quad   x \in \mathbb R,
\end{equation}
 where
\begin{eqnarray}\label{3.7}
G(x,t)=
\begin{cases}
        \ u(x)v(t), & x \geq t, \\
        \ u(t)v(x), & x \leq t, \\
     \end{cases}
\end{eqnarray}
and $u(x), ~v(x)$ are absolutely continuous positive functions defined on $\mathbb R$,  with the properties
\begin{equation}\label{3.8}
u'(t)\leq 0, \quad v'(t)\geq 0,\quad t\in \mathbb R
\end{equation}
and
\begin{equation}\label{3.9}
\lim_{t\rightarrow-\infty}v(t)=\lim_{t \rightarrow +\infty}u(t)=0.
\end{equation}

 Note that by virtue of ({\rm \ref{3.8}}) and ({\rm \ref{3.9}}), $G$ is a unimodal function  for any fixed  $x \in \mathbb R$ \cite{Cottle}.
 Therefore, it is natural to say that ({\rm \ref{3.7}}) is a uniformly unimodal function of $x\in \mathbb R$ when ({\rm \ref{3.8}) and ({\rm \ref{3.9}}) hold.
 Our goal is to find estimates of $F(x)$ for $x \in \mathbb R.$ To begin, we   introduce the functions
\begin{equation}\label{3.10}
q_{1}(t)=-\frac{u'(t)}{u(t)}, \qquad  q_{2}(t)=\frac{v'(t)}{v(t)}, \quad t \in \mathbb R.
\end{equation}
The following relations are obvious consequences of ({\rm \ref{3.8}}) and ({\rm \ref{3.9}}):
\begin{equation}\label{3.11}
 0\leq q_{1}(t),~ q_{2}(t) ~ \in L_{1}^{\textrm{loc}}(\mathbb  R),
\end{equation}
\begin{equation}\label{3.12}
 \int_{0}^{\infty}q_{1}(t)\,dt=\int_{-\infty}^{0}q_{2}(t)\,dt=\infty.
\end{equation}
Thus, we can define the functions
\begin{equation}\label{3.13}
d_{1}(x)=\inf_{d>0}\bigg\{d: \int_{x-d}^{x+d}q_{1}(t)\,dt=2 \bigg\}, \qquad x\in \mathbb{R}
\end{equation}
and
\begin{equation}\label{3.14}
d_{2}(x)=\inf_{d>0}\bigg\{d: \int_{x-d}^{x+d}q_{2}(t)\,dt=2 \bigg\}, \qquad x\in \mathbb{R}.
\end{equation}

\begin{thm}\label{t3.5} Suppose the relations  {\rm (\ref{3.11})} and {\rm(\ref{3.12})} hold, and, in addition,
 $q_{k} \in K(\gamma_{k})$,  ~ $\gamma_k \leq e^{-1}$ for $k=1,2.$ Then
\begin{equation}\label{3.15}
F(x) \asymp u(x)v(x)\left[d_{1}(x)+d_{2}(x)\right], \qquad x\in\mathbb R.
\end{equation}%
\end{thm}

An example of an application of Theorem  \ref{t3.5} is given in Section $5.$

\section{Proofs}
In what follows we will  assume that conditions (\ref{1.2}) and (\ref{1.7}) are always satisfied, and do not mention them any more.

\textbf{Proof of Theorem 3.3}

We first need to establish the following lemma.
\begin{lemma}\label{l4.1}
Suppose the equalities {\rm(\ref{2.6})}  hold for  some $a \geq  1$, $b>0$, and $x_{0}\gg 1$ such that $|x|\geq x_{0}$. Then if $q \notin L_{1}(0,\infty)$ $ \left(q \notin L_{1}(-\infty,0)\right)$,
an $R(x,bd(x))$-covering ~of $[x, \infty)$ (resp., $(-\infty, x])$  exists for any $x \in \mathbb R.$
\end{lemma}

We prove the statement of Lemma \ref{l4.1} for the $[x,~ \infty)$ semi-axis only (the case $(-\infty,~ x)$ is dealt with in a similar manner). By Theorems {\rm \ref{t2.1}} and {\rm  \ref{t2.3}}, we need  to show that ({\rm \ref{2.5}})~holds when $\varkappa(t)=bd(t)$. For $b \in (0, 1]$  we have
 $x-bd(x)\rightarrow\infty$ as $x\rightarrow\infty$, and the equality is an obvious consequence of Theorem \ref{t2.3}.
Now let $b>1$.
Let $x\geq x_{0} $ and $\{\bigtriangleup_{n}\}_{n=1}^{\infty}$ be an $R(x-bd(x), d)$-covering ~of $[x-bd(x), \infty)$.
Clearly, there are two possibilities:
    \begin{enumerate}
\item $[x-bd(x),x+bd(x)]\subseteq\bigtriangleup_{1}=[x_{1}-d(x_{1}), x_{1}+d(x_{1})];$
\item $\bigtriangleup_{1}\subset[x-bd(x),x+bd(x)].$
\end{enumerate}

 In the case 1) we have (see (\ref{1.8})):
\begin{equation}\label{4.1}
\int_{x-bd(x)}^{x+bd(x)}q(t)\,dt\leq\int_{x_{1}-d(x_{1})}^{x_{1}+d(x_{1})}q(t)\,dt=2.
\end{equation}

Consider now the case 2). It is obvious that there exists  $n>1$ such that $\bigtriangleup_{n}^{+}<x+bd(x)$ and
$\bigtriangleup_{n+1}^{+}>x+bd(x). $ For such $n$, (\ref{2.6}) yields
\begin{equation}\label{4.2}
2bd(x)\geq\sum_{k=1}^{n}2d(x_{k})=\sum_{k=1}^{n}2\frac{d(x_{k})}{d(x)}d(x)\geq\sum_{k=1}^{n}\frac{2}{a}d(x)
=\frac{2}{a}nd(x)~ \Longrightarrow ~ n\leq ab.
\end{equation}
Consequently,
\begin{equation}\label{4.3}
\int\limits_{x-bd(x)}^{x+bd(x)}q(t)\,dt\leq\sum_{k=1}^{n+1}\int_{\bigtriangleup_{k}}q(t)\,dt=2(n+1)\leq 2(ab+1).
\end{equation}
Thus, in both the cases we have:
\begin{equation}\label{4.4}
\int_{x-bd(x)}^{x+bd(x)}q(t)\,dt\leq c< \infty,\quad {\rm  for \ any ~ } |x|\geq x_{0}.
\end{equation}

Now assume that $ x-bd(x)$ does not tend  to infinity when $x \rightarrow \infty.$
Then we can find a number $c$ and a sequence $\{x_{l}\}_{l=1}^{\infty}$ such that
$x_{l}\rightarrow\infty$ as $l\rightarrow\infty,$ and $x_{l}-bd(x_{l})\leq c<\infty, ~ l=1,2,\ldots$, which in conjunction with (\ref{4.4}) implies that
\begin{equation}\notag
\infty>c \geq \int_{x_{l}-bd(x_{l})}^{x_{l}+bd(x_{l})}q(t)\,dt\geq \int_{c}^{x_{l}}q(t)\,dt~\rightarrow ~ \infty\quad  {\rm ~ as ~} \quad l\rightarrow\infty,
\end{equation}
which is a contradiction. Thus, $x-bd(x)\rightarrow\infty$ as $x \rightarrow\infty$. To finish the proof, it remains to refer to
  Theorems {\rm \ref{t2.1}} and {\rm \ref{t2.3}}.
\begin{flushright}$\blacksquare$\end{flushright}
\vspace{2mm}

\textbf{Proof of Theorem 3.4.}{}
Let us check the inequality (\ref{3.4}) for $J(x)$ (for the function $I(x)$, (\ref{3.3}) can be verified in a similar way).
In turn,  for $S(x)$  the inequality is easily deduced from a combination of the ones for $I(x)$ and $J(x)$. \\
\indent The lower estimate of $J(x)$  follows immediately from  (\ref{1.8}):
\begin{align*}
\textit{J}(x)=&\int_{x}^{\infty}\exp\bigg\{-\int_{x}^{t}q(\xi)\,d\xi\bigg\}\,dt\geq \int_{x}^{x+d(x)}\exp\bigg\{-\int_{x}^{t}q(\xi)\,d\xi\bigg\}\,dt\\&\geq
\int_{x}^{x+d(x)}\exp\bigg\{-\int_{x-d(x)}^{x+d(x)}q(\xi)\,d\xi\bigg\}\,dt
=e^{-2}d(x)=\frac{e^{-2}}{q^{*}(x)}.
\end{align*}
To  estimate $\textit{J}(x)$, $x\in \mathbb R$, from above, we distinguish three cases:
\begin{align*}
1)\quad x\geq x_{0},\qquad 2)\quad x\leq -x_{0}-1, \qquad 3)\quad x\in[-x_{0}-1, x_{0}].
\end{align*}
 The hypotheses of Theorem 3.5 are assumed to be satisfied and will not be mentioned in the following statements.
We begin with case 1)
\begin{lemma}\label{l4.3}
Let $ x\geq x_{0}.$ Then  $\textrm{J}(x)~\leq ~ c\textit{J}_{a}^{~*}(x)$, where
\begin{align}\label{4.5}
\displaystyle{\textrm{J}_{a}^{~*}(x)=\int_{x}^{\infty}\exp\bigg\{-\frac{1}{a}\int_{x}^{t}q^{*}(\xi)\,{d(\xi)}\bigg\}dt}
\end{align}
\end{lemma}
\vspace{1mm}

{\bf Proof.} Let $\{\bigtriangleup_{n}\}_{n=1}^{\infty}$ be an $R(x,d)$-covering ~of
$[x, \infty), ~ x\geq x_{0}.$ Since $b\geq1$ (see Lemma \ref{l2.7}), the inequalities
\begin{equation}\label{4.6}
a^{-1}d(x)\leq d(t)\leq ad(x), \quad t \in [x-d(x), x+d(x)],
\end{equation}%
hold for $t\geq x_{0}$.
According to Definition \ref{d2.2}, (\ref{1.8}), and (\ref{4.6}), for any $n\geq 1$ we have
\begin{align}\label{4.7}\notag
\int_{\bigtriangleup_{1}^{-}}^{\bigtriangleup_{n}^{-}}q(\xi)\,d\xi
=
&
\int_{\bigtriangleup_{1}^{-}}^{\bigtriangleup_{n}^{+}}q(\xi)\,d\xi-\int_{\bigtriangleup_{n}}q(\xi)\,d\xi
=\sum_{k=1}^{n}\int_{\bigtriangleup_{k}}q(\xi)\,d\xi-\int_{\bigtriangleup_{n}}q(\xi)\,d\xi\\&\notag=
\sum_{k=1}^{n}2-2=\sum_{k=1}^{n}\frac{1}{d(x_{k})}\int_{\bigtriangleup_{k}}1\,d\xi-2=
\sum_{k=1}^{n}\int_{\bigtriangleup_{k}}\frac{d(\xi)}{d(x_{k})}\frac{d\xi}{d(\xi)}-2\\&\geq
\sum_{k=1}^{n}\frac{1}{a}\int_{\bigtriangleup_{k}}\frac{d\xi}{d(\xi)}-2=
\frac{1}{a}
\int_{\bigtriangleup_{1}^{-}}^{\bigtriangleup_{n}^{+}}\frac{d\xi}{d(\xi)}-2.
\end{align}
Using Definition \ref{d2.2} and relating to (\ref{1.8}) and (\ref{4.7}), we get:
\begin{align*}
\textit{J}(x)
=
&
\int_{x}^{\infty}\exp\bigg\{-\int_{x}^{t}q(\xi)\,d\xi\bigg\}\,dt
=
\sum_{n=1}^{\infty}\int_{\bigtriangleup_{n}}
\exp\bigg\{-\int_{\bigtriangleup_{1}^{-}}^{t}q(\xi)\,d\xi\bigg\}\,dt
\\&
\leq
\sum_{n=1}^{\infty}2d(x_{n})\exp\bigg\{-\int_{\bigtriangleup_{1}^{-}}^{\bigtriangleup_{n}^{-}}q(\xi)\,d\xi\bigg\}
\leq
e^{2}\sum_{n=1}^{\infty}2d(x_{n})\exp\bigg\{-\frac{1}{a}\int
_{\bigtriangleup_{1}^{-}}^{\bigtriangleup_{n}^{+}}\frac{d\xi}{d(\xi)}\bigg\}
\\
&=
 c\sum_{n=1}^{\infty}\int_{\bigtriangleup_{n}}\exp\bigg\{-\frac{1}{a}\int
_{\bigtriangleup_{1}^{-}}^{t}\frac{d\xi}{d(\xi)}
-\frac{1}{a}\int_{t}^{\bigtriangleup_{n}^{+}}\frac{d\xi}{d(\xi)}\bigg\}\,dt
\leq
 c\sum_{n=1}^{\infty}\int_{\bigtriangleup_{n}}\exp\bigg\{-\frac{1}{a}\int
_{\bigtriangleup_{1}^{-}}^{t}\frac{d\xi}{d(\xi)}\bigg\}\,dt\\&=c\int_{x}^{\infty}
\exp\bigg\{-\frac{1}{a}\int_{x}^{t}\frac{d\xi}{d(\xi)}\bigg\}\,dt=c\int_{x}^{\infty}
\exp\bigg\{-\frac{1}{a}\int_{x}^{t}q^{*}(\xi)\,d\xi\bigg\}\,dt=cJ_{a}^{~*}(x).
\end{align*}%
\begin{flushright}$\blacksquare$\end{flushright}
\goodbreak

\begin{lemma}\label{l4.4}
Let $x\geq x_{0}.$  Then
\begin{equation}\label{4.8}
 \textit{J}(x)\leq cd(x)=\frac{c}{q^{*}(x)}.
 \end{equation}
\end{lemma}

{\bf Proof.}{} We use the assumptions of Lemma \ref{l4.1} to deduce that an $R(x,bd)-$
covering of $[x, \infty), ~ x\geq x_{0}$ ~ exists. Let $\{\omega_{n}\}_{n=1}^{\infty}$ be the system of segments
forming the $R(x,bd)$-covering. Then for $n\geq 1$ we have (see the notation of Lemma \ref{l4.3})
\begin{align}\label{4.9}\notag
\int_{\omega_{1}^{-}}^{\omega_{n}^{-}}\frac{d\xi}{d(\xi)}=&
\int_{\omega_{1}^{-}}^{\omega_{n}^{+}}\frac{d\xi}{d(\xi)}
-\int_{\omega_{n}}\frac{d\xi}{d(\xi)}
=\sum_{k=1}^{n}\int_{\omega_{k}}\frac{d\xi}{d(\xi)}-\int_{\omega_{n}}\frac{d\xi}{d(\xi)}
\\
&\notag=\sum_{k=1}^{n}\int_{\omega_{k}}\frac{d(x_{k})}{d(\xi)}\frac{d\xi}{d(x_{k})}-\int_{\omega_{n}}\frac{d(x_{k})}{d(\xi)}\frac{d\xi}{d(x_{k})}
\geq
 \sum_{k=1}^{n}\frac{1}{a}\int_{\omega_{k}}\frac{d\xi}{d(x_{k})}-\int_{\omega_{n}}a\frac{d\xi}{d(x_{n})}
\\&=\frac{2b}{a}n-2ab.
\end{align}
Note that the inequalities (see Definition \ref{d2.2} and Definition \ref{d2.5})
\begin{align*}
&
\frac{1}{a}\leq\frac{d(x_{n+1})}{d(\omega_{n+1}^{-})}\leq a, \quad \quad
\frac{1}{a}\leq\frac{d(\omega_{n}^{+})}{d(x_{n})}\leq a,
\end{align*}
hold for $n \geq 1$.
Since the above relations imply
\begin{align*}
\frac{1}{a^2}\leq\frac{d(x_{n+1})}{d(x_{n})}\leq a^2,
\end{align*}
we obtain
\begin{equation}\label{4.10}
\frac{1}{a^{2n-2}}\leq\frac{d(x_{n})}{d(x_{1})}\leq a^{2n-2}, \quad\quad n\geq 1.
\end{equation}
Thus, using Definition \ref{d2.2} and the inequalities (\ref{4.5}), (\ref{4.9}) and (\ref{4.10}) we get
\begin{align*}
\textit{J}(x)~\leq ~
&
 c\textit{J}_{a}^{~*}(x)= c\int_{x}^{\infty}\exp\bigg\{-\frac{1}{a}\int_{x}^{t}\frac{d\xi}{d(\xi)}\bigg\}\,dt
=
c\sum_{n=1}^{\infty}\int_{\omega_{n}}\exp\bigg\{-\frac{1}{a}\int_{\omega_{1}^{-}}^{t}\frac{d\xi}{d(\xi)}\bigg\}\,dt
\\
&\leq
c\sum_{n=1}^{\infty}d(x_{n})\exp\bigg\{-\frac{1}{a}\int_{\omega_{1}^{-}}^{\omega_{n}^{-}}\frac{d\xi}{d(\xi)}\bigg\}\leq
c\sum_{n=1}^{\infty}d(x_{n})\exp\left\{-\frac{2b}{a^{2}}n\right\}
\\
& \leq c\frac{d(x_{1})}{d(x)}\cdot d(x)\sum_{n=1}^{\infty}\frac{d(x_{n})}{d(x_{1})}\exp\left\{-\frac{2b}{a^{2}}n\right\}
\leq cd(x)\sum_{n=1}^{\infty}\gamma^{2n}=cd(x)=\frac{c}{q^{*}(x)}.
\end{align*}
\begin{flushright}$\blacksquare$\end{flushright}

Now we consider case 2).
\begin{lemma}\label{l4.4}
Let $x\leq -x_{0}.$ Then
\begin{align}\label{4.11}
 \widetilde{\textit{J}}(x)~\leq ~c~\widetilde{\textit{J}}_{a}^{~*}(x),
 \end{align}
where
\begin{equation}\label{4.12}
\widetilde{\textit{J}}(x)=\int_{x}^{-x_{0}}\exp\bigg\{-\int_{x}^{t}q(\xi)\,d\xi\bigg\}\,dt, \qquad if \quad x\leq -x_{0}
\end{equation}
\begin{equation}\label{4.13}
\widetilde{\textit{J}}_{a}^{~*}(x)=\int_{x}^{-x_{0}}\exp\bigg\{-\frac{1}{a}\int_{x}^{t}\frac{d\xi}{d(\xi)}\bigg\}\,dt, \qquad if \quad x\leq -x_{0}
\end{equation}
 \end{lemma}

{\bf Proof.} Let $\{\bigtriangleup_{n}\}_{n=-\infty}^{-1}$ be an ${\mathbb R}(-x_{0}, d)$-covering of $(-\infty, x_{0}].$
The following  estimate is obtained similarly to (\ref{4.7}), therefore we omit its derivation:
\begin{equation}\label{4.14}
\int_{\bigtriangleup_{n}^{-}}^{\bigtriangleup_{k}^{-}}q(\xi)\,d\xi~\geq~\frac{1}{a}\int_{\bigtriangleup_{n}^{-}}^{\bigtriangleup_{k}^{+}}
\frac{d\xi}{d(\xi)}-2,\qquad n\leq k\leq-1.
\end{equation}

By  the inequalities (\ref{4.12}), (\ref{4.13}), and  (\ref{4.14}),
\begin{align}\label{4.15}\notag
\widetilde{\textit{J}}(\bigtriangleup_{n}^{-})=&
\int_{\bigtriangleup_{n}^{-}}^{-x_{0}}\exp\bigg\{-\int_{\bigtriangleup_{n}^{-}}
^{t}q(\xi)\,d\xi\bigg\}\,dt=
\notag \sum_{k=n}^{-1}\int_{\bigtriangleup_{k}}\exp\bigg\{-\int_{\bigtriangleup_{n}^{-}}
^{t}q(\xi)\,d\xi\bigg\}\,dt
\\
&\notag\leq
\sum_{k=n}^{-1}2d(x_{k})\exp\bigg\{-\int_{\bigtriangleup_{n}^{-}}^{\bigtriangleup_{k}^{-}}q(\xi)\,d\xi\bigg\}\leq
e^{2}\sum_{k=n}^{-1}2d(x_{k})
\exp\bigg\{-\frac{1}{a}\int_{\bigtriangleup_{n}^{-}}^{\bigtriangleup_{k}^{+}}\frac{d\xi}{d(\xi)}
\bigg\}
\\
&
\notag= c\sum_{k=n}^{-1}\int_{\bigtriangleup_{k}}
\exp\bigg\{-\frac{1}{a}\int_{\bigtriangleup_{n}^{-}}^{t}\frac{d\xi}{d(\xi)}
-\frac{1}{a}\int_{t}^{\bigtriangleup_{k}^{+}}\frac{d\xi}{d(\xi)}\bigg\}\,dt\\&\leq
c\sum_{k=n}^{-1}\int_{\bigtriangleup_{k}}\exp\bigg\{-\frac{1}{a}\int_{\bigtriangleup_{n}^{-}}^{t}\frac{d\xi}{d(\xi)}\bigg\}\,dt=c\widetilde{\textit{J}}
_{a}^{~*}(\bigtriangleup_{n}^{-}).
\end{align}
Below we derive ({\rm \ref{4.11}}) from ({\rm \ref{4.15}}). To this end we use the  inequalities
\begin{align}\label{4.16}
\frac{2}{a}\leq\int_{\bigtriangleup(x)}\frac{d\xi}{d(\xi)}\leq2a, \qquad \bigtriangleup(x)=[x-d(x),x+d(x)], \quad x\leq -x_{0},
\end{align}
and
\begin{align}\label{4.17}
\int_{\bigtriangleup_{n}^{+}}^{t}\frac{d\xi}{d(\xi)}\geq\int_{x}^{t}\frac{d\xi}{d(\xi)}-2a, \quad x\in \bigtriangleup_{n}, \qquad t\in [\bigtriangleup_{n}^{+}, -x_{0}], \quad n\leq-1.
\end{align}
Inequality (\ref{4.16}) follows from (\ref{2.6}) (in view of Lemma \ref{l2.7}, $b\geq1$) :
\begin{align*}
\frac{2}{a}\leq\int_{\bigtriangleup(x)}\frac{d\xi}{d(\xi)}=\int_{\bigtriangleup(x)}\frac{d(x)}{d(\xi)}
\frac{d\xi}{d(x)}\leq 2a,
\end{align*}
while (\ref{4.17}) is easily obtained from (\ref{4.16}):
\begin{align*}
\int_{\bigtriangleup_{n}^{+}}^{t}\frac{d\xi}{d(\xi)}=\int_{x}^{t}\frac{d\xi}{d(\xi)}-\int_{x}^{\bigtriangleup_{n}^{+}}
\frac{d\xi}{d(\xi)}
\geq
\int_{x}^{t}\frac{d\xi}{d(\xi)}-\int_{\bigtriangleup_{n}}\frac{d\xi}{d(\xi)}
\geq
\int_{x}^{t}\frac{d\xi}{d(\xi)}-2a.
\end{align*}

Assume now that $x\in \bigtriangleup_{n}, ~ n\leq-1.$ In the next chain of calculations, we use (\ref{4.16}), (\ref{4.17}), and (\ref{4.14}):
\begin{align*}
\widetilde{\textit{J}}(x)=&\int_{x}^{-x_{0}}\exp\bigg\{-\int_{x}
^{t}q(\xi)\,d\xi\bigg\}\,dt
=
\int_{x}^{\bigtriangleup_{n}^{+}}\exp\bigg\{-\int_{x}
^{t}q(\xi)\,d\xi\bigg\}\,dt
\\&+
\int_{\bigtriangleup_{n}^{+}}^{-x_{0}}\exp\bigg\{-\int_{x}
^{t}q(\xi)\,d\xi\bigg\}\,dt\leq
 e^{2}\int_{x}^{\bigtriangleup_{n}^{+}}e^{-2}\,dt
+
\int_{\bigtriangleup_{n}^{+}}^{-x_{0}}\exp\bigg\{-\int_{\bigtriangleup_{n}^{+}}
^{t}q(\xi)\,d\xi\bigg\}\,dt
\\
&\leq
 c\int_{x}^{\bigtriangleup_{n}^{+}}\exp\bigg\{-\frac{1}{a}\int_{\bigtriangleup_{n}}\frac{d\xi}{d(\xi)}\bigg\}\,dt
+
c\int_{\bigtriangleup_{n}^{+}}^{-x_{0}}\exp\bigg\{-\frac{1}{a}\int_{\bigtriangleup_{n}^{+}}^{t}\frac{d\xi}{d(\xi)}\bigg\}\,dt
\\
&\leq
c\int_{x}^{\bigtriangleup_{n}^{+}}\exp\bigg\{-\frac{1}{a}\int_{x}^{t}\frac{d\xi}{d(\xi)}\bigg\}\,dt
+
c\int_{\bigtriangleup_{n}^{+}}^{-x_{0}}\exp\bigg\{-\frac{1}{a}\int_{x}^{t}\frac{d\xi}{d(\xi)}+2\bigg\}\,dt
\\
&\leq
c\int\limits_{x}^{-x_{0}}\exp\bigg\{-\frac{1}{a}\int_{x}^{t}\frac{d\xi}{d(\xi)}\bigg\}\,dt=
c\widetilde{\textit{J}}_{a}^{~*}(x).
\end{align*}
\begin{flushright}$\blacksquare$\end{flushright}
\begin{lemma}\label{l4.5}
Let $x\leq -x_{0}.$ Then
\begin{equation}\label{4.18}
\widetilde{\textit{J}}(x)~\leq ~cd(x)~=~ \frac{c}{q^{*}(x)}.
\end{equation}
 \end{lemma}

{\bf Proof.} Let $\{\omega_{n}\}_{n=-\infty}^{-1}$ be an $R(x, bd)$-covering of $(-\infty, -x_{0}].$
Then
\begin{align}\label{4.19}
&\int
_{\omega_{n}^{-}}^{\omega_{k}^{-}}\frac{d\xi}{d(\xi)}\geq\frac{2b}{a}|n-k|-2ab, \qquad n\leq k\leq-1,
\end{align}
\begin{align}\label{4.20}
\frac{1}{a^{2|n-k|-2}}\leq\frac{d(x_{k})}{d(x_{n})}\leq a^{2|n-k|-2}, \qquad n\leq k\leq-1.
\end{align}
This implies (similarly to Lemma \ref{l4.3})
\begin{align}\label{4.21}
\widetilde{\textit{J}}(\omega_{n}^{-})~\leq ~c \widetilde{\textit{J}}_{a}^{~*}(\omega_{n}^{-}) ~\leq ~c d(\omega_{n}^{-}), \qquad~ n\leq -1.
\end{align}
The inequalities ({\rm \ref{4.19}}), ({\rm \ref{4.20}}), ({\rm \ref{4.21}}) are checked in the same way as ({\rm \ref{4.9}}), ({\rm \ref{4.10}}), ({\rm \ref{4.8}}) respectively, therefore we drop their derivation.

Assume now that  $x\in \omega_{n}, ~ n\leq -1.$ Then, using the above relations, we have
\begin{align*}
\widetilde{\emph{J}}(x)~\leq ~ & c \widetilde{\textit{J}}_{a}^{~*}(x) = c \int_{x}^{\omega_{n}^{+}}\exp\bigg\{-\frac{1}{a}\int_{x}^{t}\frac{d\xi}{d(\xi)}\bigg\}\,dt~
\\&+
c \int_{\omega_{n}^{+}}^{-x_{0}}\exp\bigg\{-\int_{\omega_{n}^{+}}^{t}\frac{d\xi}{d(\xi)}\bigg\}\,dt\cdot
\exp\bigg\{-\frac{1}{a}\int_{x}^{\omega_{n}^{+}}\frac{d\xi}{d(\xi)}\bigg\}\\&\leq
c(\omega_{n}^{+}-x)+c\widetilde{\textit{J}}_{a}^{~*}(\omega_{n}^{+})\leq
 c\left(d(x_{n})+d(\omega_{n}^{+})\right)\leq cd(x)=\frac{c}{q^{*}(x)}.
\end{align*}

Let us continue with the case 2). To estimate  $\textit{J}(x)$ for $x\leq-x_{0}-1$, first note that the integral $\textit{J}(x_{0})$ converges.
Indeed, from the inequality $ \widetilde{\textit{J}}(x)~\leq ~cd(x)$ it follows that
\begin{align*}
\textit{J}(-x_{0})= &\int
_{-x_{0}}^{\infty}\exp\bigg\{-\int_{-x_{0}}^{t}q(\xi)\,d\xi\bigg\}\,dt
=
\int_{-x_{0}}^{x_{0}}\exp\bigg\{-\int_{-x_{0}}^{t}q(\xi)\,d\xi\bigg\}\,dt
\\&+
\int
_{x_{0}}^{\infty}
\exp\bigg\{-\int_{-x_{0}}^{t}q(\xi)\,d\xi\bigg\}\,dt\leq
2x_{0}+J(x_{0})\leq c(1+d(x_{0}))<\infty.
\end{align*}
Thus, for $x\leq -x_{0}-1$ we get
\begin{align*}
J(x)=&\int_{x}^{\infty}\exp\bigg\{-\int_{x}^{t}q(\xi)\,d\xi\bigg\}\,dt
=
\widetilde{\textit{J}}(x)+\int_{-x_{0}}^{\infty}\exp\bigg\{-\int_{x}^{t}q(\xi)\,d\xi\bigg\}\,dt
\\
&=
\widetilde{\textit{J}}(x)\Bigg\{1+\textit{J}(-x_{0})\exp\bigg\{-\int_{x}^{-x_{0}}q(\xi)\,d\xi\bigg\}\left(\widetilde{\textit{J}}(x)\right)^{-1}\Bigg\}
\\
&=
\widetilde{\textit{J}}(x)\Bigg\{1+c(x_{0})\Bigg(\exp\bigg\{\int_{x}^{-x_{0}}q(\xi)\,d\xi\bigg\}\widetilde{\textit{J}}(x)\Bigg)^{-1}\Bigg\}\\&\leq
cd(x)\Bigg\{1+\Bigg(\int_{x}^{-x_{0}}\exp\bigg\{\int_{t}^{-x_{0}}q(\xi)\,d\xi\bigg\}\,dt\Bigg)^{-1}\Bigg\}\\&\leq
cd(x)\left(1+\frac{1}{-x_{0}-x}\right)\leq cd(x)=\frac{c}{q^{*}(x)}.
\end{align*}
Hence, in both the  cases 1) and 2),  inequality (\ref{3.5}) is true.
Let us now consider the case 3).  Let
\begin{align*}
f(x)=\textit{J}(x)\cdot\frac{1}{d(x)},\quad  x\in [-x_{0}-1, x_{0}].
\end{align*}
From Theorem \ref{t2.1} it follows that  $f(x)$ is a continuous function, and,
as the interval $[-x_{0}-1, x_{0}] $ is finite and closed,
$f(x)$  is bounded.
 The statement is proved.
\begin{flushright}$\blacksquare$\end{flushright}
\vspace{2mm}

{\bf Proof of Theorem \ref{t3.1}}.
In what follows, we suppose that the assumptions of Theorem 3.1 are in force.
Let $p\in [1, \infty)$ and $f \in L_{p}$. Consider the function
\begin{align}\label{4.22}
y(x):=(Gf)(x)=\int\limits_{x}^{\infty} \exp(-\int\limits_{x}^{t}q(\xi)\,d\xi)f(t)\,dt, \quad x\in \mathbb R.
\end{align}
We claim that the integral in ({\rm \ref{4.22}}) converges for all $x\in \mathbb R$, i.e., the function $y(x)$ is defined for any $x\in \mathbb R$. Indeed, by H\"{o}lder's inequality and ({\rm \ref{3.5}}) we have
\begin{align*}
&|y(x)|\leq \int\limits_{x}^{\infty}\exp\bigg\{-\int\limits_{x}^{t}q(\xi)\,d\xi\bigg\}|f(t)|\,dt\leq
\left(\int\limits_{x}^{\infty}\exp\bigg\{-\int\limits_{x}^{t}q(\xi)\,d\xi\bigg\}\,dt\right)^{1/p'}\cdot\\&
\left(\int\limits_{x}^{\infty}\exp \bigg\{-\int\limits_{x}^{t}q(\xi)\,d\xi \bigg\}|f(t)|^{p}\,dt\right)^{1/p}\leq
c(p)d(x)^{1/p'}||f||_p<\infty, \quad x\in \mathbb R.
\end{align*}
 Obviously, $y(x)$ is an absolutely continuous function for $x\in \mathbb R$ and a solution of ({\rm \ref{1.1}}). Let us check the inclusion $y\in L_{p,\theta}$, and inequality ({\rm \ref{1.4}}). To this end, we formulate one more lemma.
\begin{lemma}\label{l4.6}
Let
\begin{equation}\label{4.23}
M(x)=\int\limits_{-\infty}^{x}\theta(t)\exp\bigg\{-\int\limits_{t}^{x}q(\xi)\,d\xi\bigg\}\,dt.
\end{equation}
Then $M=\sup_{x\in\mathbb R}M(x)<\infty$.
\end{lemma}
{\bf Proof}
Let $\{\vartriangle_{n}\}_{n=-\infty}^{-1}$ be an $R(x, bd)$-covering of $(-\infty, x]$ (see Definition 2.2 and Corollary 2.4). Then, by ({\rm \ref{1.8}}), for $n\leq -1$ we have
\begin{align}\label{4.24}
\int\limits_{\triangle_n^+}^{\triangle_{-1}^{+}}q(\xi)\,d\xi=&\int\limits_{\triangle_n^-}^{\triangle_{-1}^{+}}q(\xi)\,d\xi
-\int\limits_{\triangle_n}q(\xi)\,d\xi=\sum_{k=-n}^{-1}\int\limits_{\triangle_k}q(\xi)\,d\xi-2=\sum_{k=-n}^{-1}2-2=2(|n|-1).
\end{align}
Further, since $b\geq 1$ (see Lemma \ref{2.7}):
\begin{align}\label{4.25}
\frac{1}{a}\leq\frac{d(t)}{d(x)}\leq a \quad ~ {\rm for } ~\quad  |t-x|\leq d(x),\quad |x|\geq x_{0}.
\end{align}
But $d(x)$ is a positive continuous function for all $x\in \mathbb R$ (see Theorem \ref{t2.1}). Therefore, there exists a number $\alpha \in[1, \infty)$, such that
\begin{align}\label{4.26}
\alpha^{-1}\leq\frac{d(t)}{d(x)}\leq \alpha ~\quad {\rm for } ~ \quad  |t|\leq |x_{0}|,~ ~  |x|\leq x_{0}.
\end{align}
Denote $c=\max\{a, \alpha\}. $ It follows from ({\rm \ref{4.25}}) and ({\rm \ref{4.26}}) that
\begin{align}\label{4.27}
c^{-1}\leq\frac{d(t)}{d(x)}\leq c ~ \quad  {\rm for } ~  \quad |t-x|\leq d(x), \quad x \in \mathbb R.
\end{align}
Let us return to our Lemma. According to ({\rm \ref{4.24}}), ({\rm \ref{4.27}}), Definition 2.2, and ({\rm \ref{3.2}}), we get:
\begin{align*}
M(x)=&\int\limits_{-\infty}^{x}\theta(t)\exp\bigg\{-\int\limits_{t}^{x}q(\xi)\,d\xi\bigg\}\,dt\leq
c\int\limits_{-\infty}^{x}\frac{1}{d(t)}\exp\bigg\{-\int\limits_{t}^{x}q(\xi)\,d\xi\bigg\}\,dt\\& \notag=
c\sum\limits_{n=-\infty}^{-1}\int\limits_{\triangle_{n}}\frac{1}{d(t)}\exp\bigg\{-\int\limits_{t}^{\triangle_{-1}^{+}}q(\xi)\,d\xi\bigg\}\,dt
\\&\leq
c\sum\limits_{n=-\infty}^{-1}\bigg\{\int\limits_{\triangle_{n}}\frac{d(x_{n})}{d(\xi)}\frac{d\xi}{d(x_{n})} \bigg\}\exp\bigg\{-\int\limits_{\triangle_{n}^{+}}^{\triangle_{-1}^{+}}q(\xi)\,d\xi\bigg\} \\&\notag \leq c\sum\limits_{n=-\infty}^{-1}\exp\bigg\{-2(|n|-1)\bigg\}=c<\infty
\end{align*}
\begin{flushright}$\blacksquare$\end{flushright}
Below we obtain inequality ({\rm \ref{1.4}}). Let $f\in L_{p}, ~p\in[1, \infty)$. Using  H\"{o}lder's inequality, ({\rm\ref{3.5}}), ({\rm \ref{3.2}), Lemma \ref{l4.6}, and Fubini's theorem we get
\begin{align*}
&||y(x)||_{p,\theta}^p=\int\limits_{-\infty}^{\infty}\theta(x)^p\left\vert\int\limits_{x}^{\infty}\exp\bigg\{-\int\limits_{x}^{t}q(\xi)\,d\xi\bigg\}|f(t)|^p\,dt\right\vert^p\,dx
\\
&\notag\leq
\int\limits_{-\infty}^{\infty}\theta(x)^p\left(\int\limits_{x}^{\infty}\exp\bigg\{-\int\limits_{x}^{t}q(\xi)\,d\xi\bigg\}\,dt\right)^{p/p'}
\left(\int\limits_{x}^{\infty}\exp\bigg\{-\int\limits_{x}^{t}q(\xi)\,d\xi\bigg\}|f(t)|^{p}\,dt\right)\,dx
\\
&\notag\leq
c\int\limits_{-\infty}^{\infty}\theta(x)^pd(x)^{p-1}\left(\int\limits_{x}^{\infty}\exp\bigg\{-\int\limits_{x}^{t}q(\xi)\,d\xi\bigg\}|f(t)|^p\,dt\right)\,dx
\\&\notag\leq
c\int\limits_{-\infty}^{\infty}\theta(x)\left(\int\limits_{x}^{\infty}\exp\bigg\{-\int\limits_{x}^{t}q(\xi)\,d\xi\bigg\}|f(t)|^p\,dt\right)\,dx= c\int\limits_{-\infty}^{\infty}M(t)|f(t)|^p\,dt\leq c||f||_{p}^{p}.
\end{align*}
Hence, it remains to check that a solution of the homogenous equation
\begin{align}\label{4.28}
-z'(x)+q(x)z(x)=0, \quad x \in \mathbb R
\end{align}
belongs to the space $L_{p,\theta}$ only if $z \equiv 0$. Clearly, we can write the solution of (\ref{4.28}) as

\begin{align}\label{4.29}
z(x)=\alpha \exp\bigg\{\int\limits_{x_{0}}^{x}q(\xi)\,d\xi\bigg\}, \quad x \in \mathbb R
\end{align}
where $\alpha$ is constant and $x_{0}$ is as in Definition \ref{2.5}. Let $\alpha \neq 0$ and , with this, $z \in L_{p,\theta}$. The following estimates are based on (4.29), (2.6), (3.2) , (4.7), (4.9); also, $R[x_{0} , d)$ and $R[x_{0}, bd)-$ coverings of $[x_{0}, \infty)$ are used.
\begin{align}\label{4.30}\notag
&\infty>||z(x)||_{p,\theta}^p\geq\int\limits_{-\infty}^{\infty}|\theta(t)z(t)|^p\,dt\geq
(|\alpha|c^{-1})^p\int\limits_{x_{0}}^{\infty}\frac{1}{d(t)^p}\exp\bigg\{\int\limits_{x_{0}}^{t}q(\xi)\,d\xi\bigg\}\,dt\\&\notag=
(|\alpha|c^{-1})^p\sum\limits_{n=1}^{\infty}\int\limits_{\triangle_n}\frac{1}{d(t)^p}\exp\bigg\{p\int\limits_{\triangle_{1}^{-}}^{t}q(\xi)\,d\xi\bigg\}\,dt\geq
(|\alpha|c^{-1})^p\sum\limits_{n=1}^{\infty}\left(\int\limits_{\triangle_n}\frac{dt}{d(t)^p}\right)\exp\bigg\{p \int\limits_{\triangle_{1}^{-}}^{\triangle_{n}^{-}}q(\xi)\,d\xi\bigg\}\\&\notag\geq
(|\alpha|c^{-1})^p\sum\limits_{n=1}^{\infty}\left(\int\limits_{\triangle_n}\frac{dt}{d(t)^p}\right)\exp\bigg\{\frac{p}{a}\int\limits_{\triangle_{1}^{-}}^{\triangle_{n}^{+}}\frac{d\xi}{d(\xi)}\bigg\}
\geq(|\alpha|c^{-1})^p\sum\limits_{n=1}^{\infty}\int\limits_{\triangle_n}\frac{1}{d(t)^p}\exp\bigg\{\frac{p}{a}\int\limits_{\triangle_{1}^{-}}^{t}\frac{d\xi}{d(\xi)}\bigg\}\,dt
\\&\notag=
(|\alpha|c^{-1})^p\int\limits_{x_0}^{\infty}\frac{1}{{d(t)}^{p}}\exp\bigg\{\frac{p}{a}\int\limits_{x_0}^{t}\frac{d\xi}{d(\xi)}\bigg\}\,dt
=
(|\alpha|c^{-1})^p\sum\limits_{n=1}^{\infty}\int\limits_{\omega_n}\frac{1}{d(t)^p}\exp\bigg\{\frac{p}{a}\int\limits_{\omega_{1}^{-}}^{t}\frac{d\xi}{d(\xi)}\bigg\}
\,dt\\&\notag\geq
(|\alpha|c^{-1})^p\sum\limits_{n=1}^{\infty}\left(\int\limits_{\omega_n}\left(\frac{d(x_n)}{d(t)}\right)^p\frac{dt}{{d(x_n)}^p}\right)\exp\bigg\{\frac{p}{a}\int\limits_{\omega_{1}^{-}}^{\omega_{n}^{-}}\frac{d\xi}{d(\xi)}\bigg\}
\geq(|\alpha|c^{-1})^p\sum\limits_{n=1}^{\infty}\frac{1}{d(x_n)^{p-1}}\exp\bigg\{2\frac{b}{a^2}pn\bigg\}\\& := (|\alpha|c^{-1})^p\sum\limits_{n=1}^{\infty}
u_n.
\end{align}
Using d'Alambert criterion for the convergence of series (\ref{4.30}) we get (see \ref{4.10})
\begin{align*}
\frac{u_{n+1}}{u_n}=\left(\frac{d(x_{n})}{d(x_{n+1})}\right)^{p-1}\exp\bigg\{2\frac{b}{a^2}p\bigg\}\geq\frac{a^2}{a^{2p}}\exp\bigg\{\frac{2b}{a^2}p\bigg\}\geq(\frac{1}{\gamma})^{2p}\geq e^{2p}>1.
\end{align*}
This is a contradiction. Hence, $\alpha=0$ and $z\equiv 0$. The statement is proved.
\begin{flushright}$\blacksquare$\end{flushright}
{\bf Proof of Theorem 3.5}
The following relations are  obvious consequences of (\ref{3.10}), (\ref{3.11}), (\ref{3.12}), and Theorem 3.3:
\begin{align*}
F(x)=
&
\int_{-\infty}^{\infty}G(x,t)\,dt=u(x)\int_{-\infty}^{x}v(t)\,dt+v(x)\int_{x}^{\infty}u(t)\,dt
\\
&=
 u(x)v(x)\left(\int_{-\infty}^{x}\frac{v(t)}{v(x)}\,dt+\int_{x}^{\infty}\frac{u(t)}{u(x)}\,dt\right)
\\
&=
u(x)v(x)\left(\int_{-\infty}^{x}\exp\bigg\{-\int_{t}^{x}q_{2}(\xi)\,d\xi\bigg\}\,dt+\int_{x}^{\infty}\exp\bigg\{-\int_{x}^{t}q_{1}(\xi)\,d\xi\bigg\}\,dt\right)
\\
&
\asymp u(x)v(x)\left(d_{1}(x)+d_{2}(x)\right), \quad x \in \mathbb R.
\end{align*}
\begin{flushright}$\blacksquare$\end{flushright}
{\bf Proof of Theorem 1.3}{}
{\it Necessity.} Suppose equation (\ref{1.1}) is correctly solvable in $L_p$ for $p \in [1, \infty)$. Then, according to Theorem 1.2 , $q_{0}(a)>0$ for some $a \in (0, \infty)$, and (\ref{1.7}) holds. Let $n$ be the smallest integer such that $nq_{0}(a)\geq 1$. Then for any $x \in \mathbb R$ we have
\begin{align*}
\int\limits_{x-2na}^{x+2na}q(\xi)\,\xi=\int\limits_{x-2na}^{x}q(\xi)\,d\xi+\int\limits_{x}^{x+2na}q(\xi)\,d\xi\geq
nq_{0}(a)+nq_{0})a)\geq 2\ \quad \Rightarrow
\end{align*}
\begin{align*}
 d(x)\leq 2na \quad \Rightarrow \quad d_0 < \infty
\end{align*}
{\it Sufficiency.} Since (\ref{1.7}) holds, $d(x)$ is defined for any  $x \in \mathbb R$, and
\begin{align*}
\int\limits_{x-d_0}^{x+d_0}q(\xi)\,d\xi \geq \int\limits_{x-d(x)}^{x+d(x)}q(\xi)\,d\xi=2.
\end{align*}Hence, (\ref{1.5}) holds.
\begin{flushright}$\blacksquare$\end{flushright}

\section{Examples}
The following examples are applications of Theorems 3.1 and 3.5.

{\it Example 1.}{}
Consider the equation
\begin{align}\label{5.1}
-y'(x)+q(x)y(x)=f(x), \quad x \in \mathbb  R,
\end{align}
where $f \in L_p$,~ $p\in[1, \infty)$,
\begin{align}\label{5.2}
q(x)=\frac{1}{(1+x^2)^{\alpha}}+\frac{\cos(1+x^2)^\beta}{(1+x^2)^\alpha}, \quad x \in \mathbb R,
\end{align}
and the parameters $\alpha, ~ \beta$ satisfy the conditions
\begin{align}\label{5.3}
0<\alpha<\frac{1}{2}, \quad \alpha+ \beta> \frac{1}{2}.
\end{align}
We formulate our results in the following theorem.
\begin{thm}\label{t5.1} Let $\theta(x)$ be a positive continuous function for $x \in \mathbb R$ which satisfies the relation
$$\theta(x)\asymp (1+x^2)^{-\alpha}, \quad x \in \mathbb R. $$
 Then the pair $\{L_{p, \theta}; L_p\}$ is admissible  for equation {\rm ( 5.1)} for any $p\in [1, \infty)$.
\end{thm}

{\bf Proof.}{} Clearly, the requirements (\ref{1.2}), (\ref{1.6}) hold. It remains to check (\ref{3.1}). Since (\ref{5.2}) is an even function, we investigate it on the $[0, \infty)$ half-axis only. We introduce the sequences
$$x_{k}=\sqrt{(2\pi k)^{1/\beta}-1} \quad{\rm and } \qquad z_{k}=\sqrt{(2\pi k+\frac{\pi}{4})^{1/\beta}-1}, ~\qquad k \geq 1.$$
Let $\delta=\frac{1}{2}-\alpha$ and $\varkappa_{0} \gg 1$. In the following estimates we use the binomial series:
 \begin{align*}
 &\int\limits_{0}^{\infty}q(t)\,dt\geq\sum\limits_{k=k_0}^{\infty}\int\limits_{x_k}^{z_k}\left(
 \frac{1}{(1+x^2)^\alpha}+\frac{\cos(1+x^2)^\beta}{(1+x^2)^\alpha}\right)\,dx\\&\geq
 \sum\limits_{k=k_0}^{\infty}\int\limits_{x_k}^{z_k}
 \frac{dx}{(1+x^2)^\alpha}\geq \sum\limits_{k=k_0}^{\infty}\frac{z_k-x_k}{(1+z_k^2)^\alpha}=
 \sum\limits_{k=k_0}^{\infty}\frac{z_k^2-x_k^2}{(1+z_k^2)^\alpha}\frac{1}{z_k+x_k}
 \\&\geq c^{-1}\sum\limits_{k=k_0}^{\infty}\frac{(2\pi k+\frac{\pi}{4})^{{1}/{\beta}}-(2\pi k)^{1/{\beta}}}
 {(1+z_k^2)^\alpha k^{1/2\beta}}\geq
 c^{-1}\sum\limits_{k=k_0}^{\infty}\frac{ k^{1/\beta}\left((1+\frac{1}{8k})^{{1}/{\beta}}-1\right)}{(1+z_k^2)^\alpha \cdot k^{1/2\beta}}\\&\geq
 c^{-1}\sum\limits_{k=k_0}^{\infty}\frac{ 1}{(1+z_k^2)^\alpha}\frac{1}{k^{1-1/2\beta}}\geq
  c^{-1}\sum\limits_{k=k_0}^{\infty}\frac{ 1}{k^{(\alpha/\beta)+1-1/2\beta}} \\&\notag= c^{-1}\sum\limits_{k=k_0}^{\infty}\frac{ 1}{k^{1-\delta/\beta}}=\infty.
   \end{align*}
   Thus,  (\ref{3.1}) holds, and $d(x)$ is defined for all $x \in \mathbb R$ (see Theorem \ref{t2.1}).
   Let us now find the exact order two-sided estimates of $d(x)$. To this end we use Theorem 2.8. Let   \begin{align*}
q_{1}(x)=\frac{1}{(1+x^{2})^\alpha} \quad  q_{2}(x)=\frac{\cos(1+x^{2})^{\beta}}{(1+x^{2})^{\alpha}}, \quad~x\in \mathbb{R},
\end{align*}
in (\ref{2.8}).  In view of the obvious relations
 \begin{align}\label{5.4}
[x-2(1+x^{2})^\alpha, ~x+2(1+x^{2})^\alpha]\subseteq \sigma(x)=[x-4x^{2\alpha},~ x+4x^{2\alpha}]
\end{align}
 \begin{align}\label{5.5}
\lim\limits_{x\rightarrow-\infty}(x+4x^{2\alpha})=-\infty,\quad \lim\limits_{x\rightarrow\infty}(x-4x^{2\alpha})=\infty,
\end{align}
which hold for $|x|\geq x_0\gg1$, we have
\begin{align}\label{5.6}
\frac{1}{2}\leq 1-\frac{4}{|x|^{1-2\alpha}}\leq|\frac{s}{x}|\leq 1+\frac{4}{|x|^{1-2\alpha}}, \quad s \in \sigma(x),
\end{align}and
\begin{align}\label{5.7}
\frac{1}{16}\leq \frac{1}{4}\left(\frac{s}{x}\right)^2\leq \frac{1+s^2}{1+x^2}\leq 4\left(\frac{s}{x}\right)^2\leq 16, \quad s \in \sigma(x).
\end{align}
Below to estimate $\varkappa_1(x)$ when $|x|\rightarrow \infty$ we use relations (\ref{5.6}) and (\ref{5.7}):
\begin{align*}
\varkappa_{1}(x)=&(1+x^{2})^{2\alpha}\sup\limits_{|\xi|\leq 2(1+x^{2})^{\alpha}} \left\vert\int\limits_{x-\xi}^{x+\xi}\left(\frac{1}{(1+s^{2})^{\alpha}}\right)''\,ds\right\vert
\leq c|x|^{4\alpha}\sup\limits_{|\xi|\leq 4|x|^{2\alpha}} \left\vert\int\limits_{x-\xi}^{x+\xi}\frac{ds}{(1+s^{2})^{\alpha+1}}\right\vert\\&\notag\leq
c\frac{|x|^{6\alpha}}{(1+x^{2})^{\alpha+1}}\leq\frac{c}{|x|^{2-4\alpha}}\rightarrow 0 ~{\rm ~when }~~ |x|\rightarrow \infty.
\end{align*}
Let us now estimate $\varkappa_2(x)$ for $|x|\gg 1$. Note that
$f(s)=s(1+s^{2})^{\beta+\alpha-1}$ is a monotonically increasing function when $s\gg  1$ (see 5.3). Hence,
using (5.4), (5.5), (5.6), (5.7) and the second mean value theorem (see \cite{Titch}), we get:
\begin{align*}
\varkappa_{2}(x)=&\sup\limits_{|\xi|\leq 2(1+x^{2})^{\alpha}} \left\vert\int\limits_{x-\xi}^{x+\xi} \frac{\cos(1+s^{2})^\beta}{(1+s^{2})^{\alpha}}\,ds\right\vert \leq
\sup\limits_{|\xi|\leq 2(1+x^{2})^{\alpha}} \left\vert\int\limits_{x-\xi}^{x+\xi} \frac{(\sin(1+s^{2})^\beta)'}{2\beta s(1+s^{2})^{\alpha+\beta -1}}\,ds\right\vert\\&\leq
\frac{c}{|x|(1+x^{2})^{\beta+\alpha-1}}\sup\limits_{|\xi|\leq 2(1+x^{2})^{\alpha}} \left\vert\sin{(1+s^{2})^{\beta}}|_{x-\xi}^{x+\xi}\right|\vert \leq
\frac{c}{|x|^{2(\beta+\alpha)-1}}
\rightarrow 0\quad{\rm ~as}\quad |x|\rightarrow \infty.
\end{align*}
Thus, by Theorem \ref{t2.7}, we have
\begin{align*}
c^{-1}(1+x^2)^\alpha\leq d(x)\leq c(1+x^2)^\alpha, \quad x \in \mathbb R,
\end{align*}and
\begin{align*}
d(x)=(1+\varepsilon(x))(1+x^2)^\alpha, \quad |\varepsilon(x)|\leq \frac{c}{|x|^\nu}, \quad |x|\gg 1.
\end{align*}
Here $\nu=\min\{4\alpha, 2\beta+2\alpha-1\}$.
Further, conditions (\ref{2.13}) and (\ref{2.14}) are obviously satisfied. So by Theorem \ref{t2.8},
$q \in K(\gamma)$ , $ \gamma \leq e^{-1}$. Now Theorem \ref{t5.1} is seen to be a consequence of Theorem \ref{t3.1}, and the proof is finished.
\begin{flushright}$\blacksquare$\end{flushright}

{\it Example 2.}{}
Consider the integral
$F_{1}(x)=\int\limits_{-\infty}^{\infty} G_{1}(x,t)\,dt, ~  x \in \mathbb{R},$ \label{5.1} where
\begin{eqnarray}\label{5.8}
G_{1}(x,t)=
\begin{cases}
        \ \exp\left\{t^{3}+t\cos{t} -x^{3}\right\}, & t \leq x \\
        \ \exp\left\{x^{3}- t\cos{t} -t^{3}\right\}, & t \geq x. \\
     \end{cases}
\end{eqnarray}
Our goal is to find  the exact order two-sided estimates of $F_{1}(x)$ for all $x\in \mathbb R.$
The question is reduced to estimation of the integral of the type (\ref{3.4}). The following relations are obvious:
 \begin{align*}
F_{1}(x)&=\int_{-\infty}^{x}G_{1}(x,t)\,dt+\int_{x}^{\infty}G_{1}(x,t)\,dt\\
&=e^{x\cos{x}}\int_{-\infty}^{x}e^{-(x^{3}-t^{3}+x\cos{x}-t\cos{t})}\,dt\\
&+e^{-x\cos{x}}\int_{x}^{\infty} e^{-(t^{3}-x^{3}+t\cos{t}-x\cos{x})}\,dt\\
&=e^{x\cos{x}}\int_{-\infty}^{x} e^{-\int_{t}^{x}(3\xi^{2}+(\xi\cos{\xi})')\,d\xi}\,dt\\
&+e^{-x\cos{x}}\int_{x}^{\infty} e^{-\int_{x}^{t}(3\xi^{2}+(\xi\cos{\xi})')\,d\xi}\,dt\\
&=e^{x\cos{x}}\int_{-\infty}^{x}e^{-\int_{t}^{x}(3\xi^{2}-\xi\sin{\xi})\,d\xi}\cdot e^{(\sin{t}-\sin{x})}\,dt\\
&+e^{-x\cos{x}}\int_{x}^{\infty} e^{-\int_{x}^{t}(3\xi^{2}-\xi\sin{\xi})\,d\xi} \,\cdot e^{(\sin{x}-\sin{t})}\,dt.
\end{align*}
Hence,
\begin{align*}
e^{-2}h(x,t)< F_{1}(x)<e^{2}h(x,t),
\end{align*}
where
\begin{align}\label{5.9}
h(x,t)= e^{x\cos{x}}\int_{-\infty}^{x} e^{-\int_{t}^{x}(3\xi^{2}-\xi\sin{\xi})\,d\xi}\,dt+
 e^{-x\cos{x}}\int _{x}^{\infty} e^{-\int_{x}^{t}(3\xi^{2}-\xi\sin{\xi})\,d\xi} \,dt.
\end{align}

It remains to estimate the integrals
 \begin{align}\label{5.10}
 I_1(x)=\int _{-\infty}^{x}e^{-\int_{t}^{x}q(\xi)\,d\xi}\,dt   \quad{\rm and}\quad J_1(x)=\int_{x}^{\infty}e^{-\int_{x}^{t}q(\xi)\,d\xi}\,dt
 \end{align}
 for all $x\in \mathbb{R},$  where $q(t)=3t^{2}-t\sin{t}$, $t\in \mathbb R$.
In order to apply Theorem \ref{t3.3}, let us introduce the functions (see Theorem \ref{t2.7})
\begin{align}\label{5.11}
q_{1}=3t^{2}+1,\quad q_{2}=-1-t\sin{t},\quad t\in \mathbb{R},
\end{align}
and check that $q \in K(\gamma), ~ \gamma\leq e^{-1}$.
  Now for $\varkappa_1(x),~ x \gg 1$, we have
\begin{align*}
\varkappa_{1}(x)=\frac{1}{q_{1}(x)^{2}}\sup_{|\xi|\leq 2/q_{1}(x)} \left\vert\int_{x-\xi}^{x+\xi}q_{1}''(s)\,ds\right\vert
\leq \frac{c}{|x|^{6}}\rightarrow 0
\end{align*}
as $|x|\rightarrow \infty.$
To estimate $\varkappa_2(x)$ when $|x| \rightarrow  \infty$, we use the second mean value theorem (see \cite{Titch}):
\begin{align*}
\varkappa_{2}(x)&=\sup_{|\xi|\leq 2/q_{1}(x)} \left\vert\int_{x-\xi}^{x+\xi}q_{2}(s)\,ds\right\vert =
\sup_{|\xi|\leq 2/3x^{2}+1} \left\vert\int_{x-\xi}^{x+\xi}(1+t\sin{t})\,dt\right\vert \\&\notag\leq
\frac{c}{x^{2}}+\sup_{|\xi|\leq 2/3x^{2}+1}\left\vert\int\limits_{x-\xi}^{x+\xi}t\sin{t}\,dt\right\vert\leq
\frac{c}{x^{2}}+\sup_{|\xi|\leq 2/3x^{2}+1} \sup_{\theta \in [x-\xi, x+\xi]}|x+\xi|\left\vert\int\limits_{\theta}^{x+\xi}\sin{t}\,dt\right\vert\\&\notag\leq
\frac{c}{x^2}+c|x|
\sup_{|\xi|\leq 2/3x^{2}+1} \sup_{\theta \in [x-\xi, x+\xi]}\left\vert\cos{(x+\xi)}-\cos{\theta}\right\vert\\&\notag\leq
\frac{c}{x^2}+c|x|
\sup_{|\xi|\leq 2/3x^{2}+1} \sup_{\theta \in [x-\xi, x+\xi]}\left\vert\sin{\frac{x+\xi-\theta}{2}}\right\vert
\\&\notag\leq
\frac{c}{x^{2}}+c\frac{|x|}{3x^{2}+1}\rightarrow 0, \qquad {\rm as}  \qquad |x|\rightarrow \infty.
\end{align*}
Hence,
$$d(x)=\frac{1+\varepsilon(x)}{3x^2+1}, \qquad
|\varepsilon(x)|\leq\frac{c}{|x|},\quad |x|\gg 1,$$
so
$$d(x)\asymp \frac{1}{3x^2+1},\quad x\in \mathbb{R}.$$
The hypotheses of Theorem \ref{t2.8}  are obviously satisfied, hence, $q \in K(\gamma), ~ \gamma \leq e^{-1}$. Finally, using Theorem \ref{t3.3} we get
\begin{align*}
I_{1}\asymp J_{1}\asymp \frac{1}{3x^{2}+1} \asymp \frac{1}{x^{2}+1}, \quad  x \in \mathbb{R}
\end{align*}
and
\begin{align*}
F_{1}(x)\asymp  \frac{\cosh(x\cos{x})}{x^{2}+1}, \quad x \in \mathbb{R}.
\end{align*}

\end{document}